\documentclass[12pt]{article}

\usepackage{empheq}
\usepackage[numbered,framed]{matlab-prettifier}
\usepackage{url}
\usepackage{enumitem}
\usepackage{listings}
\usepackage{adjustbox}
\usepackage[margin=1in,footskip=0.25in]{geometry}
\usepackage{relsize}
\usepackage{amsthm}%
\usepackage{mathtools}
\usepackage{amsfonts}%
\usepackage{pdfpages}
\usepackage{amsmath}
\usepackage{color}
\usepackage{flexisym}

\definecolor{mygreen}{RGB}{28,172,0} 
\definecolor{mylilas}{RGB}{170,55,241}

\usepackage{cite}
\newtheorem{theorem}{Theorem}
\newtheorem{remark}[theorem]{Remark}
\newtheorem{conj}[theorem]{Conjecture}
\newtheorem{definition}[theorem]{Definition}

\newtheorem{proposition}[theorem]{Proposition}

\begin{document}

\title{Chaotic hidden attractor in a fractional order system modelling the interaction between dark matter and dark energy}

\author{Marius-F. Danca{\footnote{Corresponding author}}\\
STAR-UBB Institute, Babes-Bolyai University\\
400084 Cluj-Napoca, Romania\\
Email: m.f.danca@gmail.com
}

\maketitle

\begin{abstract}
In this paper the dynamics of a fractional order system modelling the interaction between dark matter and dark energy is analytically and numerically studied. It is shown for the first time that systems modelling the interaction between dark matter and dark energy, chaotic hidden attractors can be present. The chaotic attractor coexists with two asymptotically stable equilibria. Equilibria of the linearized system exhibit a center-like behavior. The numerical integration is done by means of the Adams-Bashforth-Moulton scheme and the finite Lyapunov exponents are numerically determined with a dedicated Matlab code. The 3D representation of the chaotic hidden attractor reveals the fact it is not connected with the equilibria, being ``hidden'' somewhere in the considered phase space.
\end{abstract}

\textbf{keyword }Fractional calculus; Caputo fractional derivative; Hidden chaotic attractor

\vspace{3mm}

\section{Introduction}
The interaction between dark matter and dark energy, which affects the cosmic structures, is one of not completely solved problems in cosmology and has been studied extensively \cite{unux}. Our universe, which continues to expand, consists not only of matter but also dark matter and dark energy. Today, we know that our universe is composed of about 68.3 dark energy, 26.8 dark matter and 4.9 regular matter according to Planck data \cite{doix} (detailed references on interaction between dark energy and dark matter can be found in \cite{unux} and \cite{doix}). In \cite{doix} is introduced a coupled system of Integer Order (IO) modeling the Interaction between Dark Matter and Dark Energy (IDMDE) (see also \cite{ayd2} where the coupling interactions among dark energy, dark matter, and matter are modeled by three coupled logistic models)

\begin{equation}\label{eq1}
\begin{aligned}
\dot{x}_1=&x_2x_3-x_1,\\
\dot{x}_2=&(x_3-p)x_1-x_2,\\
\dot{x}_3=&1-x_1x_2,
\end{aligned}
\end{equation}
where $p$ is a real parameter considered positive in this paper. Despite his simple form, in \cite{doix} it is shown that the system has a clear physical sense and the numerical study realized reveals that for $p=3.46$, the system behaves chaotically.

Without being one of the algebraically simplest chaotic systems (see \cite{spr1,spr2}), its rich potential chaotic behavior represented the motivation to go beyond the obtained results in \cite{doix} and to consider the Fractional Order (FO) form of the IO system \eqref{eq1}, study which reveal more interesting dynamics such as the existence of hidden attractors which to our knowledge it's a novelty and, therefore, could influence the theory of interaction coupling between dark matter and dark energy and help to understand the observed acceleration of cosmic expansion.

The concept of non-integer or FO derivatives dates back to the beginning of the theory of differential calculus, while the development of fractional calculus dates since the work of Euler, Liouville, Riemann, Letnikov \cite{treix,patrux}. The fractional calculus is as old as the integer-order one, but with application exclusively in mathematics.
 Fractional calculus can be applied to real life problems and, therefore, fractional derivatives and integrals are useful in engineering and mathematics, being helpful for scientists and researchers working with real-life applications (see, e.g., \cite{doix,cincix}). One of the famous operator of FO is Caputo's fractional derivative, proposed by Caputo in 1967 and utilized in this paper. Fractional derivatives model phenomena which takes account of interactions within the past (equation with ``memory''). Moreover, fractional derivatives describe nonlocality in space and time. Basic aspects of the theory of fractional
differential equations can be found in the monograph \cite{sasex}, or in \cite{treispex}. Also see \cite{saptex} for a review of definitions of fractional derivatives and other operators. One of the usually utilized numerical schemes for integration of continuous systems of FO which is used in this paper is the Adams-Bashforth-Moulton predictor-corrector method \cite{doispex}.

Not only fractional differential equations aroused the interest of scientists, but also fractional difference equations, one of the first definition of a fractional difference operator dating from 1974 \cite {optx} (several references on fractional differences can be found e.g. in \cite{dadus}). Initial value problems for fractional differences are described in \cite{optopt}, while applications can be found in \cite{patru,nouanoua,zecezece,unspeunspe}.

Both for continuous time and discrete time systems, Caputo's differential operator allows the formulation of initial conditions for initial-value problems in traditional form, as for integer-order initial-value problems, reason why it was chosen in this paper.

It is to note that the dynamics of many real systems are better described non-integer derivatives which incorporate memory and hereditary features of systems. So, in order to accurately describe these systems, the fractional-order differential equations have been introduced.
In chaos theory, it has been observed that chaos occurs in dynamical systems of order 3 or more. With the introduction of fractional-order systems, it was shown that chaos in systems FO systems exist even the total order is less than 3.

Recently a new term \emph{hidden attractor} has been introduced in IO systems and, later, in FO systems (see e.g. \cite{zecex} or \cite{nouax}, one of the first works on this subject). If the
attraction basin of a considered attractor intersects with any open neighborhoods of an equilibrium, the attractor is called \emph{self-excited}, otherwise, it is called \emph{hidden}. Usually, the self-excited attractors are numerically derived from unstable equilibria, while hidden attractors are difficult to be localized, because their attraction basins have no relation with small neighborhoods of any equilibria. References on systems type where hidden attractors appear are presented in \cite{dadi1}. The importance of hidden attractors is underlined for example by the hidden oscillations which might occur at aircrafts and launchers control systems \cite{unspex}.

The paper focus on showing numerically that a FO variant of an existing IO system modelling the interaction between dark matter and dark energy, can present hidden attractors. For this purpose, it is proved numerically and via Matignon's criterion, that equilibria are asymptotically stable for all fractional order $q\in(0,1)$ and parameter $q$. Next, based on the existence of local Lyapunov Exponents and their three-dimensionally plot, it is shown that the attracting set, reached by trajectories starting from an appropriate basin, is chaotic. The LEs are determined with a dedicated Matlab code. Since all considered trajectories attracted by the chaotic attractor do not cross relatively small neighborhoods of the stable equilibria, the attractor verifies the definition of a hidden attractor. Also, it is shown that, considering too short time integration intervals (approach used in several works), could lead to false results, i.e., chaotic transients, like in the considered system. A FO variant of divergence for the considered system is presented as an interesting approach of stability.

The  paper is organized as follows: Section 2 introduces the FO variant of the IDMDE system; Section 3 presents equilibria of the FO system and their stability; Section 4 deals with Lyapunov exponents; In Section 5 the chaotic hidden attractor is analyzed and the Conclusion section ends the paper.

\section{IDMDE system of FO}

We consider in this paper the commensurate FO variant of the IDMDE system \eqref{eq1} in Caputo's sense, modeled by the following FO initial value problem
\begin{equation}\label{eqx1}
D_*^qx(t)=f(x(t)),~~~ x(0)=x_0,
\end{equation}
where $x=(x_1,x_2,x_3)$, and $f:\mathbb{R}^3\rightarrow \mathbb{R}^3$ is the vector-valued function
\[f(x)=
\begin{pmatrix}
   x_2x_3-x_1\\
    (x_3-p)x_1-x_2\\
    1-x_1x_2
\end{pmatrix},
\]
and $D_*^q$ stands as Caputo's derivative with starting point $0$, of order $q$ being defined as follows
\begin{equation*}\label{cap}
D_*^qx=\frac{1}{\Gamma \left \lceil{q}\right \rceil-q }\int_0^t(t-\tau)^{\left \lceil{q}\right \rceil -q-1}D^{\left \lceil{q}\right \rceil }x(\tau)d\tau.
\end{equation*}
For $\left \lceil{q}\right \rceil\in \mathbb{N}$, $D^{\left \lceil{q}\right \rceil}$ becomes the standard differential operator of IO and for $q\in(0,1)$, as considered in this paper, $D^{\left \lceil{q}\right \rceil}x=x\textprime$.

With the above notations, the FO IDMDE system becomes

\begin{equation}\label{eq2}
\begin{aligned}
D^q_*{x}_1=&x_2x_3-x_1,\\
D^q_*{x}_2=&(x_3-p)x_1-x_2, ~~~ x(0)=x_0,\\
D^q_*{x}_3=&1-x_1x_2,
\end{aligned}
\end{equation}
with $q\in(0,1)$ and $p\in \mathbb{R}^+$.

The deeper numerical and theoretical analysis of the FO system \eqref{eq2} made in this paper reveals more interesting dynamics than those presented in \cite{doix}, for example: the pure complex conjugated eigenvalues of equilibria indicate the use of center manifold approach, while attractors coexistence phenomenon leads to the existence of chaotic hidden attractors.

\section{Equilibria}
Equilibria, the two real solutions of the parametric equation $f(x)=0$, $p\in \mathbb{R}^+$, are:
\[
X_{1,2}^*=\left(\mp\frac{\sqrt{2}}{2}\sqrt{P+p},\pm\frac{\sqrt{2}}{4}\Bigl( p-P\Bigr)\sqrt{P+p},\frac{1}{2}\left(p+P\right)     \right)
\]
where
\[
P:=\sqrt{p^2+4}.
\]
The Jacobian is
\[J(x)=
\begin{bmatrix}

    -1&x_3&x_2\\
x_3-p&     -1&x_1\\
-x_2&-x_1&0
\end{bmatrix}.
\]

Without restricting the generality, consider $p=5$. Then, equilibria are
\[
X_{1,2}^*(\mp2.2787,\mp0.4388,5.1926),
\]
and, for both equilibria, the eigenvalues are
\[\label{sd}
\sigma_1=- 2, ~~\sigma_{2,3}=\pm 2.32053i.
\]
The fact that $\operatorname{Re}(\sigma_{2,3})=0$, suggests that $X_{1,2}^*$ should be analyzed via center manifold theory. However, Matignon Theorem for autonomous FO systems \cite{lulu} can still be applied to determine the stability of $X_{1,2}^*$.
Therefore, the following result can be proved
\begin{theorem} \label{th1}For $p=5$, equilibria $X_{1,2}$ are asymptotically stable for $q\in(0,1)$.
\begin{proof}
 Denote with
 \begin{equation}\label{ioio}
 \iota =q-2\frac{|\alpha_{min}|}{\pi},
 \end{equation}
where $\alpha_{min}$ stands as the minimum of the arguments of eigenvalues. As stated by the stability theorem \cite{lulu}, the considered equilibrium is asymptotically stable if $\iota<0$, otherwise it is unstable.
For $p=5$, the arguments of the eigenvalues are $\alpha_1=\pi$, $\alpha_{2,3}=\pm\pi/2$ and $|\alpha_{min}|=\pi/2$. Then
\[\label{io}
 \iota =q-1<0,~~\text {for}~q\in(0,1),
 \]
 and, therefore, equilibria $X_{1,2}$ are asymptotically stable.
 \end{proof}
 \end{theorem}

 Because, for the general case of $p\in[p_{min},p_{max}]$, with $p_{min,max}>0$, the explicit form of eigenvalues is too complicated to be presented here, let the following conjecture sustained by computation evidence according to which the stability result Theorem \ref{th1} is true for all considered values of $p$.
 \begin{conj}\label{conj1}
 Equilibria $X_{1,2}^*$ are asymptotically stable for $q\in(0,1)$ and all $p>0$.
\end{conj}
\noindent\emph{Numerical proof.}

\begin{figure}
\begin{center}
\includegraphics[scale=0.75]{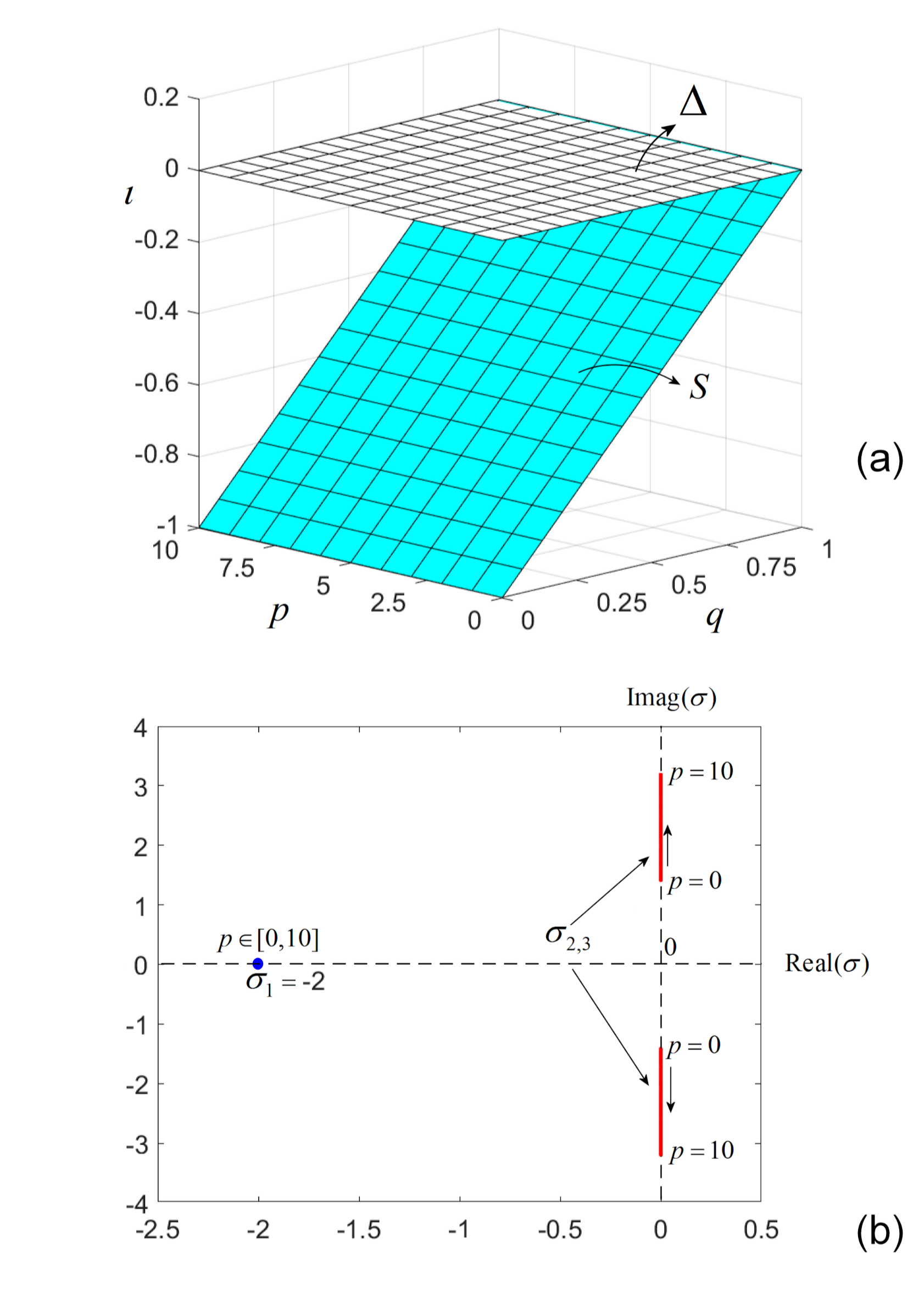}
\caption{Stability of equilibria $X_{1,2}^*$. (a) Graph $S$ generated by $\iota$ \eqref{ioio} as function of $p$ and $q$ (light blue); (b) Eigenvalues plot as function of $p$: the real eigenvalue $\sigma_1=-2$ (blue plot), and the pure complex conjugate eigenvalues $\sigma_{2,3}$ (red plot) for $p\in[0,10]$.}
\label{fig1}
\end{center}
\end{figure}

Consider lattice $L=\{(p,q)|p\in [p_{min},p_{max}], q\in(0,1)\}$ with $p_{min}=0$ and $p_{max}=10$. Then, exploring $L$ with a finite relative small step, following the reasoning in Theorem \ref{th1} and using the symbolic calculus to determine the spectrum of the eigenvalues at every grid point $(p,q)\in L$, one can obtain $\iota$ with relation \eqref{ioio} which determines a two parametric surface in the space $(p,q,\iota)$, denoted by $S $ (Fig. \ref{fig1} (a)). Note that for every $q$ one eigenvalue is real constant: $\sigma_1=-2$, while the other two are complex conjugate with $\operatorname{Re}(\sigma_{2,3})=0$ (Fig. \ref{fig1} (b)). As can be seen in Fig. \ref{fig1} (a), the graph of $S$ is situated under the horizontal critical plane $\Delta$ ($\iota=0$), for $q\in(0,1)$ and all values of $p$, which means that all points of $S$ have coordinate $\iota$ negative. Moreover, from Fig. \ref{fig1} (a) one deduces that for every fixed value of $q\in(0,1)$, $\iota$ does not depend on $p\in[p_{min},p_{max}]$.
 Therefore, for $q\in(0,1)$, $X_{1,2}^*$ are stable for all considered values of $p$.
\qed

 Note that the surface $S$ crosses the plane $\iota=0$ for $q=1$ (case not considered in this paper), the stability of $X_{1,2}^*$ being critical.

\begin{proposition}\label{propilus}
For all considered $p$, the system \eqref{eqx1} is dissipative.
\begin{proof}
The generalized divergence $\nabla(f(x))=\mathrm{Tr} (J(x))=-2<0$ and not depending on $p$. Therefore, the system is dissipative for all $p$.
\end{proof}
\end{proposition}
To note that the dissipativeness could be useful to study the rapport between the absorbed and supplied energy of the system.
 \begin{remark}
 \begin{itemize}
 \item [i)] Since $Real(\sigma_{2,3})=0$, in the IO variant of the system \eqref{eqx1}, the stability of equilibria $X_{1,2}^*$ cannot be studied via the known negativeness criterion of real components of eigenvalues and the stability has to be analyzed via center manifold theory. However, it is difficult to get the analytical expression of the centre manifold, but there are methods of the construction or approximation of center manifold (see e.g. \cite{wiki}). A useful web service which constructs center manifolds for autonomous IO systems with general linearisation is \cite{wiki2};
 \item [ii)] The existence of center manifolds and its approximation for autonomous FO systems defined with Caputo's derivative, is studied in \cite{center} and the center manifold of the Lorenz system of FO is calculated in \cite{center2}.
     \item[ii)] Related to Conjecture \ref{conj1} and Proposition \ref{propilus}, it is worth noting that for IO systems, in \cite{rus1,rus2} the instability condition of a point $x$, $\nabla f (x) > 0$, is considered (see also \cite{rus3}).  Let consider here the FO variant of divergence, which in our case, when $f=(f_1,f_2,f_3)^T$ and $x=(x_1,x_2,x_3)$, has the form \cite{jj}
     \[
     \nabla^q f(x)=D_{*x_1}^qf_1(x_1,x_2,x_3)+D_{*x_2}^qf_2(x_1,x_2,x_3)+D_{*x_3}^qf_3(x_1,x_2,x_3),
          \]
     where $D_{*{x}}^q:=\frac{\partial^q }{\partial x^q}$ represents the partial Caputo's derivative with respect $x$.  $D_*^q$ being a linear operator, if $D_*^q f(x)$ and $D_*^q g(x)$ exist, then for $a,b\in \mathbb {R}$ one has
     \[
     D_*^q(af(x)+bg(x))=aD_*^q f(x)+bD_*^q g(x),
     \]
     and also
     \[
     D_*^qx^n=\frac{\Gamma(n+1)}{\Gamma(n-q+1)}x^{n-q}, ~~~\text{for}~~ n\in \mathbb{Z}^+~~~ \text{and~~} D_*^q const=0.
     \]
     Therefore, for $q\in(0,1)$ and $p\in[p_{min},p_{max}]$, within a rectangular planar domain $(x_1,x_2)$ which contains the coordinates $x_1^*$ and $x_2^*$ of equilibria $X_{1,2}^*$ (since $f_3$ does not depend on $x_3$, $\frac{\partial^q }{\partial x_3^q }f_3(x_1,x_2)=0$, and $\nabla^qf(x)$ also does not depend on $x_3$), one has
     \[
     \nabla^qf(x)=-\frac{\Gamma(2)}{\Gamma(2-q)}(x_1^{1-q}+x_2^{1-q})<0, ~~\text{for }~~q\in(0,1).~~~ 
     \]
     In Fig.\ref{fig2}, $\nabla^qf(x)$ is determined as function of $x_1$ and $x_2$, for $q=0.995$.
 \end{itemize}
 \end{remark}

\section{Local Lyapunov exponents}

While for a given system ``global'' Lyapunov Exponents (LEs) offer information about the average growth of small perturbations, finite-time (or local) LEs provide ``local'' growth rates along a finite-time section of the trajectory and reveal the separation between trajectories. In this paper, finite-time LEs are considered.

Consider the general FO initial value problem \eqref{eqx1} with $f:\mathbb{R}^n\rightarrow \mathbb{R}^n$ a smooth function. Then, there exists the time-variant variational system associated to \eqref{eqx1} is the following linear FO system \cite{exista}

\begin{equation}\label{eqqq}
\begin{aligned}
&D^q_*\Phi(t)=J(t)\Phi (t),\\
&\Phi(0)=I_n,
\end{aligned}
\end{equation}
where $\Phi\in \mathbb{{R}}^{n\times n}$ is the fundamental solution to \eqref{eqqq} and $J(t)$ is the $n\times n$ Jacobian of $f$ evaluated at the point $x(t)$, solution of \eqref{eqx1}, i.e. along the trajectory of the system \eqref{eqx1}.

The following existence result can be enounced (the proof can be found in \cite{le1})

\begin{theorem}\label{the}
The LEs are defined by the eigenvalues of the fundamental solution $\Phi$ obtained by solving the following combined (extended) system
\begin{equation}\label{eqq}
\left\{
\begin{array}{rcr}
    D_*^q x(t) \\
    D_*^q \Phi(t)
  \end{array}
\right\}
   =\left\{
 \begin{array}{ccr}
  f(x(t)) \\
  J(t)\cdot\Phi(t)
  \end{array}
 \right\}, ~~
   \left \{
  \begin{array}{ccr}
    x(0)\\
    \Phi(0)
    \end{array}
   \right\}=
   \left \{
  \begin{array}{ll}
    x_0 \\
    I_n
  \end{array}
   \right\}.
\end{equation}
\end{theorem}

The explicit form of \eqref{eqq} allows the numerical determination of LEs by using some numerical scheme for FDEs. Simultaneously with the integration of the system \eqref{eq2} (first equation of the extended system \eqref{eqq}), one integrate the second (variational) system in \eqref{eqq}, where the solution of the first system is used.
To obtain LEs, in this paper the Matlab code FO\_Lyapunov for commensurate FO systems \cite{le22,le1} is utilized (for non-commensurate case see \cite{le2,le223}). The code adapts the Benettin procedure \cite{bene} based on Multiplicative Ergodic Theorem of Oseledec \cite{tradi} and using the Gram-Schmidt ortogonalization procedure to FO systems \cite{le1}. Despite the low efficacy and high sensitivity on algorithm parameters, and the potential loss of ortogonality, which all might lead to inaccurate numerical results, the algorithm proposed by Benettin remains one of most used numerical algorithms for IO LEs and seems to be useful for FO LEs too.

For the numerical integration consider the discretization of the time-integration interval $I=[0,T]$, $T>0$, on which the numerical solution is determined, with grid points of some equidistant partition of $I$, $t_i=hi$, $i=0,1,2,...,N$, where $h$ is a fixed step size, $h=T/N$. With these ingredients in this manuscript one considers the known predictor-corector method Adams-Bashforth-Moulton (ABM) for FO systems \cite{doispex}\footnote{A fast Matlab implementation can be found at \cite{roby}.}.

In this paper the LEs along of a chaotic trajectory are determined (for the stable attractors $X_{1,2}^*$, obviously, LEs are negative for all $q$ and $p$). For this purpose, the maximum time integration interval is $[0,T]$, with $T=1700$, the integration step size $h=0.02$ and the one of the most important parameter, $h_{norm}$, in FO\_Lyapunov representing the time moments when the Gram-Schmidt ortogonalization is applied \cite{le1}, was chosen multiple of $h$. Compared to several other cases of FO systems where the FO\_Lyapunov code has been applied, in this case the value of $h_{norm}$ is relatively higher, but chosen together with $h$ so as to be consistent: positive or negative values of LEs correspond, in the vast majority of cases, to chaotic or stable attractors, respectively.

As known, longer time integration intervals in the integration of IO systems could lead to wrong results (see e.g. \cite{lung1,lung2}). However, on the other side, this system presents chaotic transients sets which must be avoided (see e.g. the case $p=1.2$ and $q=0.995$ (point $M$ in Fig. \ref{fig3} (d)) for the initial condition $x_0=(0.1,0.1,0.1)$, with the corresponding time series in Fig. \ref{fig3} (e) where, after a chaotic transient which lasts at about $t=1300$, the trajectory tends to the first coordinate of the stable equilibrium $X_1^*(-1.3290,-0.7525,1.7662)$). Therefore, the intensive realized numerical experiments proved that the relative large time interval $[0,1700]$, represents the best choice to overcome chaotic transients.

Due to these mentioned characteristics of the FO IDMDE system and to inherent numerical errors which characterizes the Benettin algorithm, amplified by the numerical errors introduced by the ABM method for FO systems, the results should be interpreted with caution. For this system and presented data, the observed errors of LEs are of order of $5e-2$.

Denote by $\lambda_1$, $\lambda_2$, $\lambda_3$ with $|\lambda_1|<|\lambda_2|<|\lambda_3|$ the spectrum of the LEs. Because LEs depend on $q$ and also on $p$, they can be plotted in the space $(p,q,\lambda)$ as function of $p$ and $q$ (see Fig. \ref{fig3} (a), where the most representative range for $q$, $q\in[0.98,1)$ has been considered; for $q\in(0,0.98)$ nothing dynamically interesting happens).
To every $\lambda_i$, $i=1,2,3$ it corresponds a surface denoted $S_\lambda$. In this way one can obtain a general view of LEs, not only as the usual representation as function of the parameter $p$, or as function of $q$. The surface $S_{\lambda_2}$ (blue plot) is under the horizontal plane $\lambda=0$ and close to it indicating the negativeness of $\lambda_2$ for all considered $q$ and $p$. Similarly, $S_{\lambda_1}$ (green plot) represents the negative values of $\lambda_1$. The only surface presenting positive values (red plot) is $S_{\lambda_3}$.

As expected, the LEs along trajectories leading to equilibria $X_{1,2}^*$ are all negative and are not considered here.
\begin{remark}\label{rem6}
In the combined system \eqref{eqq}, $x(t)$ stands as the current state of the trajectory of the underlying system \eqref{eq2}, used to find $\Phi$ and, therefore, the initial condition $x_0$ should be chosen inside the basin of attraction of the studied attractor (here chaotic) otherwise, wrong results might be obtained.
As known, if $p$ or $q$ are varied while \eqref{eqq} is integrated, the attraction basin might change and, therefore, the reference trajectory calculated by the algorithm for LEs could``jump'' to another trajectory leading to wrong LEs  (see the switched between positive (corresponding to the chaotic attractor) and negative or zero values (corresponding to the stable equilibria $X_{1,2}^*$) in Fig. \ref{fig3} (b) and (c) which, by the mentioned reasons, couldn't be avoided). For the simulations utilized in this paper to obtain the LEs, the utilized initial conditions are $x_0=(1e-3,1e-3,1e-3)$ (note that other initial conditions could enter other attraction basins, see Fig. \ref{fig3} (e) where the trajectory is obtained with the initial condition $x_0=(0.1,0.1,0.1)$).
    \end{remark}

\section{Chaotic hidden attractor}

 Considering that the positiveness of at least one LE ensures the existence of chaos\footnote{It is possible to construct examples in which the system has positive Lyapunov exponents along a zero solution of the original system but, at the same time, this zero solution of original nonlinear system is Lyapunov stable (Perron's effect\cite{pero}, see also \cite{wikus}).}, one can see that the chaotic attractor, exists only for a relative small domain in the plane $(p,q)$ (see the surface $S_{\lambda_3}$ (red plot), situated above the plane $\lambda=0$ in Fig. \ref{fig3} (a) and also the projection of this positive part of $S_{\lambda_3}$ on the plane $(p,q)$, denoted $B_{HA}$, in Fig. \ref{fig3} (d), which can be considered as an attraction parametric basin-like). Beside the parameters $(p,q)$ which generate chaotic dynamics, $B_{HA}$ also allows the determination of minimum value of $q$ for which chaos is enhanced, namely a value nearly $0.988$.

Summarising the results on LEs (Section 4), from Figs. \ref{fig3} (a) and (d), the following important property of the FO system \eqref{eq2} can be enounced
\begin{proposition}\label{propo}
For all considered values of $p\in [0,10]$, there exists $q\in(0,1)$ for which the system \eqref{eqx1} behaves chaotic. Moreover, for $q$ close to 1, the system is chaotic for all $p$.
\end{proposition}

\begin{figure}
\begin{center}
\includegraphics[scale=0.6]{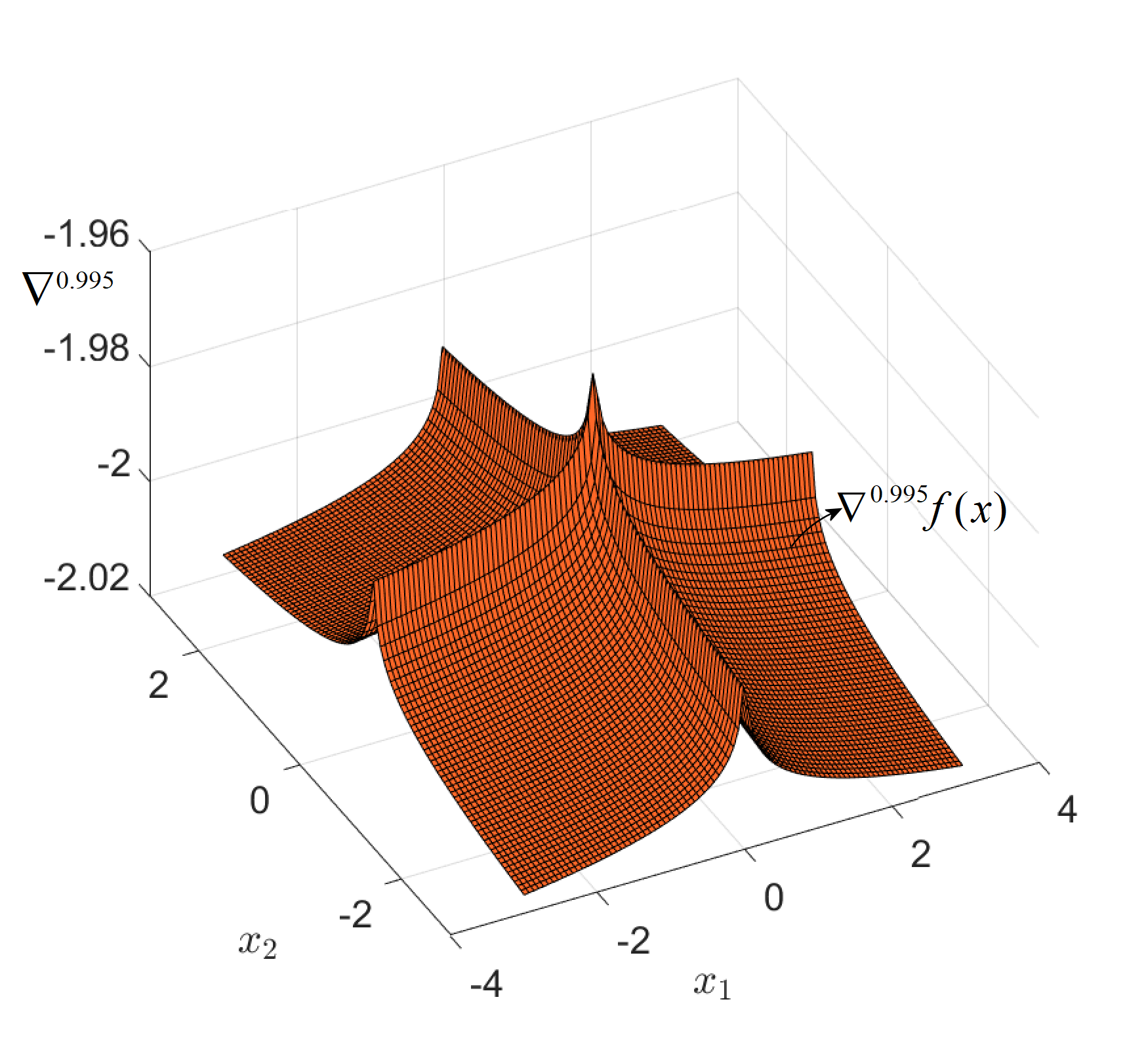}
\caption{Graph of the fractional divergence $\nabla^{0.995} f(x)$.}
\label{fig2}
\end{center}
\end{figure}

\begin{figure}
\begin{center}
\includegraphics[scale=0.6]{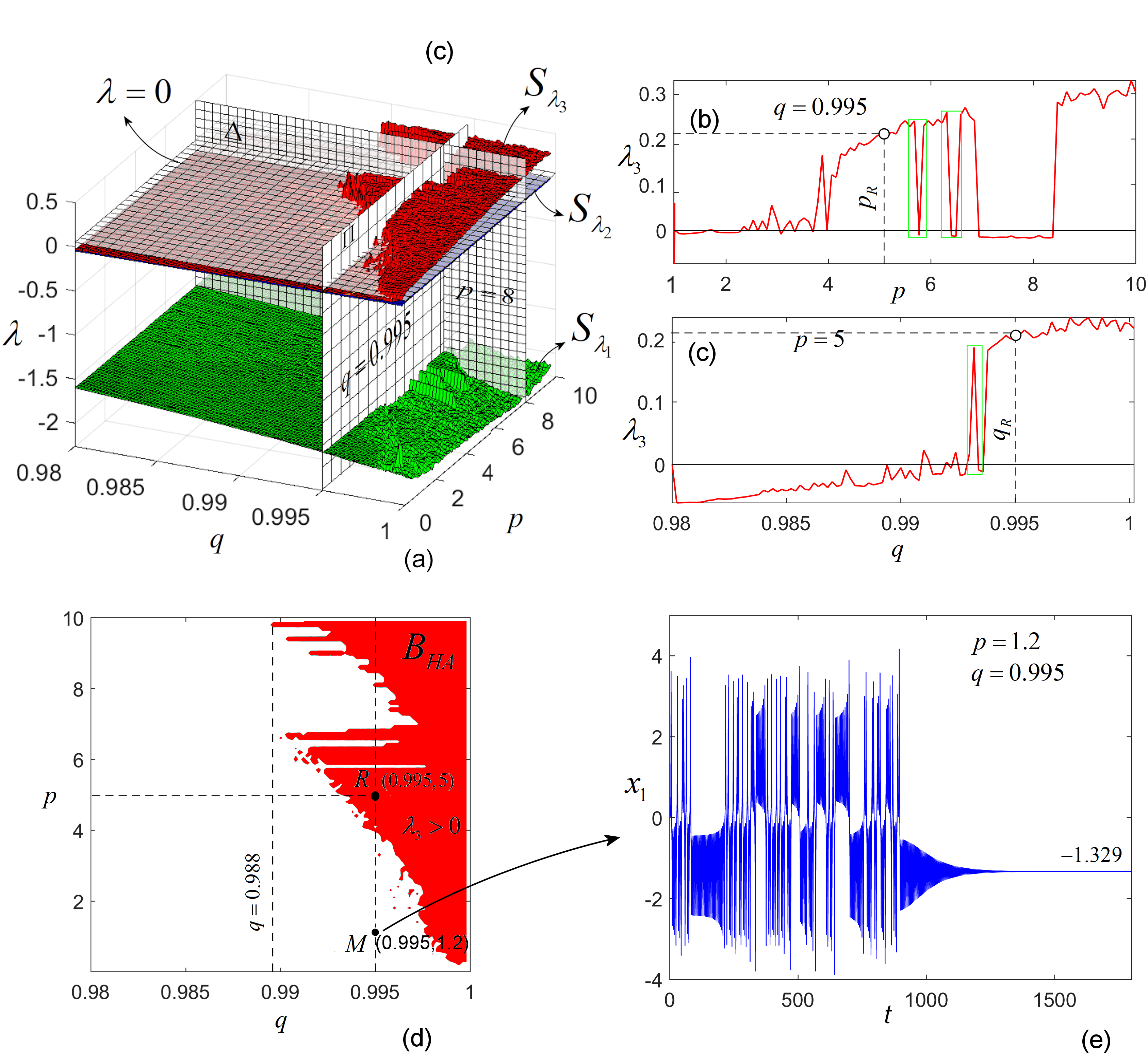}
\caption{Local LEs. (a) Surfaces $S_{\lambda_i}$, $i=1,2,3$, determined by LEs, plotted as function of $p$ and $q$; (b) Evolution of the maximum LE, $\lambda_3$, as function of $p$ for $q=0.995$ (see Fig. \ref{fig3} (d) for $p_R$);  (c) Graph of the maximum LE, $\lambda_3$, as function of $q$ for $p=5$  (see Fig. \ref{fig3} (d) for $q_R$); (d) Projection on the plane $(p,q)$ of the positive surface $S_{\lambda_3}$ (red plot), denoted $B_{HA}$; (e) Time series of a relative longue chaotic transient corresponding to $p=1.2$, $q=0.995$ corresponding to point $M$ in Fig. \ref{fig3} (d) for the initial condition $(0.1,0.1,0.1)$.}
\label{fig3}
\end{center}
\end{figure}

If one intersect the surfaces $S_{\lambda_i}$, $i=1,2,3$, with an orthogonal plane to the plane $(p,q)$ and parallel with the plane $(p,\lambda)$, the intersection curves indicate the evolution of $\lambda_i$, $i=1,2,3$ as function on $p$ (see plane $\Pi$ through $q=0.995$ in Fig. \ref{fig3} (a)). Similarly, one can obtain the evolution of LEs, for some fixed value $p$, as function of $q$, if one intersect an orthogonal plane parallel to the plane $(q,\lambda)$ (see plane $\Delta$ through $p=8$ in Fig. \ref{fig3} (a)). In Fig \ref{fig2} (b) the variation of the maximal LE, $\lambda_3$, as function of $p$, determined for $q=0.995$, $\lambda_3$, is consistent with the intersection of $S_{\lambda_3}$ with the plane $\Pi$. Fig. \ref{fig3} (c) represents the evolution of the maximal LE, $\lambda_3$, for $p=5$, as function of $q$.

\begin{definition}\cite{nouax,zecex}\label{defi}
An attractor is called a hidden attractor if its basin of attraction does not intersect with a certain open neighbourhood of equilibrium points; otherwise it is called a self-excited attractor.
\end{definition}
Following the stability of equilibria $X_{1,2}^*$ (Conjecture \ref{conj1}), and also due the existence of an chaotic attractor (Proposition \ref{propo}), the following result can be established

\begin{figure}
\begin{center}
\includegraphics[scale=0.65]{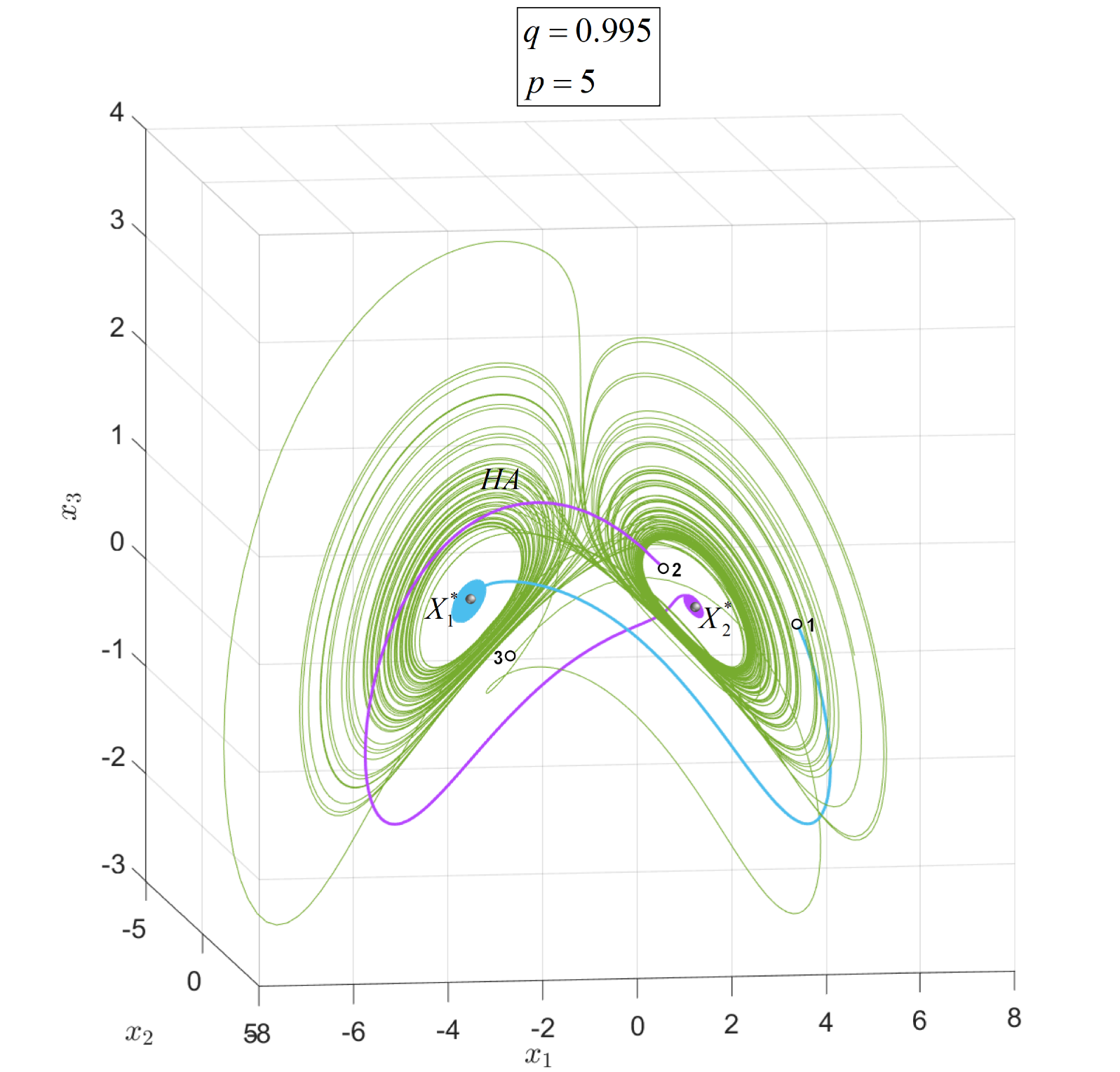}
\caption{The chaotic hidden attractor, HA (dark green plot), and the stable equilibria $X_{1,2}^*$ (magenta plot) for $p=5$ and $q=0.995$.  }
\label{fig4}
\end{center}
\end{figure}

\begin{proposition}\label{propi}
The chaotic attractor of the system \eqref{eqx1} is hidden.
\begin{proof}
Following Definition \ref{defi}, the proof is obvious since $X_{1,2}^*$ are asymptotically stable and attract trajectories starting from neighborhoods of $X_{1,2}^*$.
\end{proof}
\end{proposition}
In Fig.\ref{fig4} is presented the chaotic hidden attractor (dark green plot), denoted HA, corresponding to $p=5$ and $q=0.995$.

\begin{figure}
\begin{center}
\includegraphics[scale=0.6]{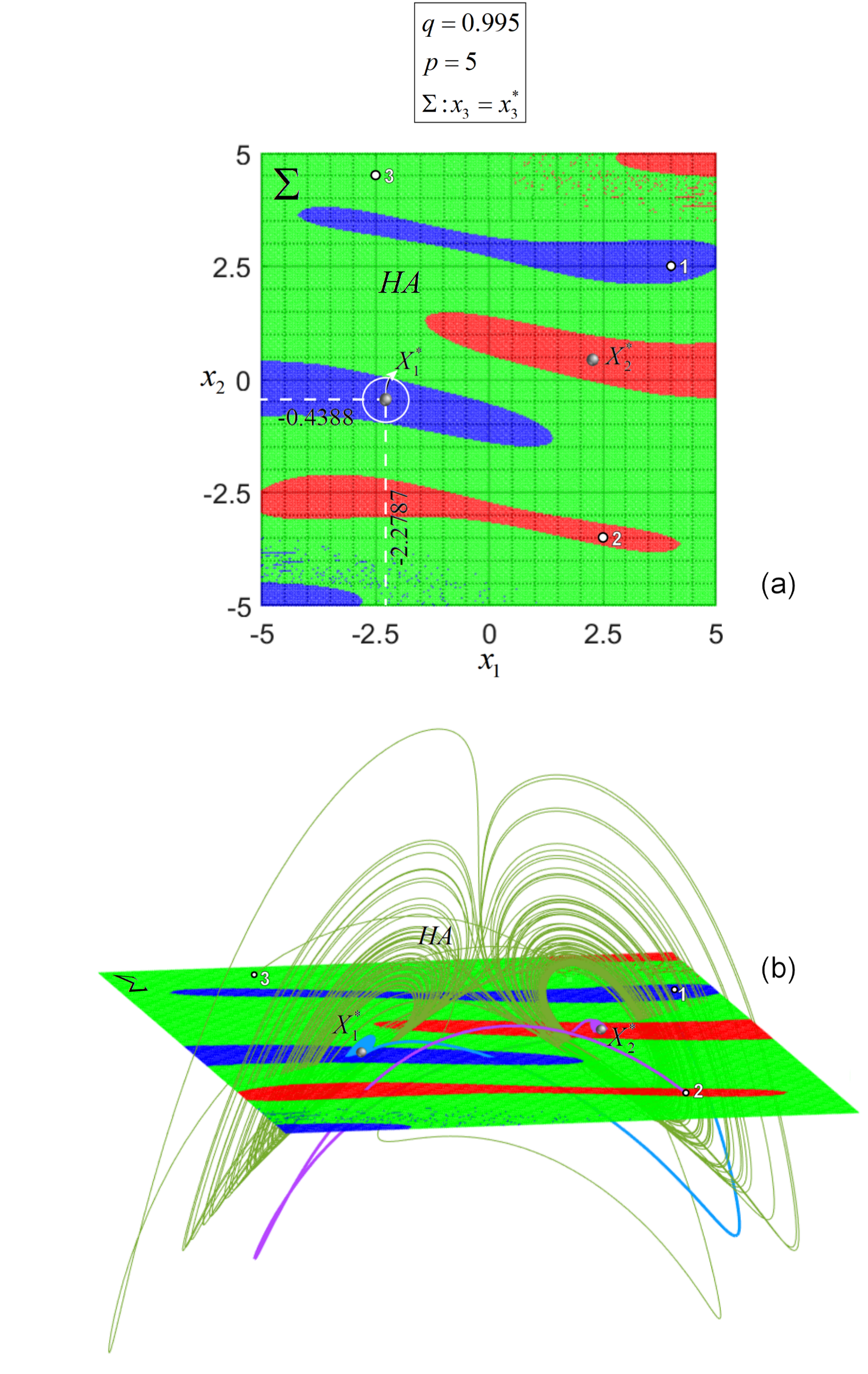}
\caption{Attraction basin of the chaotic HA for $p=5$ and $q=0.995$. (a) Horizontal intersection $\Sigma$ of the attraction basin of the chaotic HA, with the plane $x_3=x_3^*$, containing equilibria $X_{1,2}^*$. Blue points are attracted by the equilibrium $X_1^*$, red points by the equilibrium $X_2^*$, while green points represent the initial conditions in the plane $x_3=x_3^*$ leading to HA; (b)  The trajectory generating the chaotic HA (dark green plot)  crosses $\Sigma$ only on light green points. Points denoted 1,2 and 3 represent the initial conditions $x_0$ used to generate trajectories tending to $X_1^*$, $X_2^*$ and HA, respectively (see also Fig. \ref{fig5} (a)). }
\label{fig5}
\end{center}
\end{figure}

To sustain numerically Proposition \ref{propi}, let a two-dimensional representation of the attraction basin of HA corresponding to $p=5$ and $q=0.995$, consider the plane $\Sigma=\{(x_1,x_2,x_3^*)|x_1,x_2\in \mathbb{R}, x_3^* =5.1926\}$ which contains equilibria $X_{1,2}^*$ (Fig. \ref{fig5} (a)). Because third coordinates of $X_{1,2}^*$, $x_3^*$, are identic $\Sigma$ is horizontal and passes through $X_{1,2}^*$. In this plane a lattice of $100\times100$ points $(x_1,x_2,x_3^*)\in[-5,5]\times[-5,5]\times \{x_3^*\}$, is considered and at each of these points, considered as initial condition $x_0$, the system \eqref{eqx1} is integrated to determine where trajectories tend. Points which lead to equilibria $X_{1,2}^*$ are plotted blue (corresponding to $X_1^*$) or red (corresponding to $X_2^*$), while points which tend to the chaotic attractor are plotted light green. As can be seen, around equilibria there exists a relative large neighborhoods (withe disc in Fig. \ref{fig5} (a)) which contain only points attracted by one of equilibria $X_{1,2}$ respectively\footnote{A three-dimensional representation of an attraction basin of chaotic hidden attractors is presented in \cite{graf}.}.

HA, plotted in Fig. \ref{fig4}, starts from the initial condition $x_0=(-2.5,4.5,x_3^*)$, $x_3^*=5.1926$ (point $3$ in Figs \ref{fig5}), but following the definition of an attractor, could be generated from any other ``green'' initial condition. The stable attractors $X_{1,2}^*$ (light blue and magenta plots), respectively, are reached from initial condition $x_0=(4,2.5,x_3^*)$ (point $1$ in Figs \ref{fig5}) and initial condition $x_0=(2.5,-3.5,x_3^*)$ (point $2$ in Fig. \ref{fig5}). Similarly to the chaotic attractor HA, for the considered case, horizontal section through $x_3^*$, $X_{1,2}^*$ can be reached from any other ``blue'' or ``red'' initial condition.

 In the parametric space, HA corresponds to the point $R$ in Fig. \ref{fig2} (d), having the corresponding maximum LE, $\lambda_3$, marked by the black circles in Figs. \ref{fig3} (b) and (c). As expected, for these values of $p$ and $q$, the value of $\lambda_3$ is similar in both LE representations, $\lambda_3\approx 2.1$.

As expected, Fig. \ref{fig5} (b) shows that the chaotic HA crosses the horizontal plane through $x_3=x_3^*$ only in light green points, which represent the attraction basin of HA, and not red or blue points which represent the attraction basins of $X_{1,2}^*$.

As can be seen in the crossing plane $\Sigma$ (Fig. \ref{fig5} (b)) the area of the attraction basin of $HA$ is quite larger than the area of attraction basins of stable equilibria $X_{1,2}^*$, fact which suggests that the probability of the system to evolve chaotically for arbitrary initial conditions is quite large.

\section*{Conclusion}
In this paper a novel IO system modelling the interaction between dark matter and dark energy has been considered in the general form of a fractional order system modeled by Caputo's derivative. Its rich dynamics are estimated theoretically and numerically. Because the interaction between dark matter and dark energy affecting the cosmic structures is one of not completely solved problems in cosmology, the fractional approach as generalization of dynamics of integer order, could help the understanding of these dynamics such as the presence of chaos for total order less than 3, or the existence of hidden attractors.
The system presents a center equilibria whose stability could be analyzed via Matignon criterion, not via the sign study of the real components of related eigenvalues. The numerical integration has be done with the ABM method for FO systems and the local Lyapunov exponents are found with the Matlab code FO\_Lyapunov.m. Because the equilibria are asymptotical stable, the existing chaotic attractor is hidden. Further analysis of the FO system for larger values of the parameter $p$ has to be realized.

\section*{Acknowledgement} Author thanks Abdullah Gokyildirim and Haris Calgan for verifying the existence of the chaotic HA by means of circuit implementations.

\newpage{\pagestyle{empty}\cleardoublepage}

\newpage{\pagestyle{empty}\cleardoublepage}

\newpage{\pagestyle{empty}\cleardoublepage}

\end{document}